\documentclass{article}
\usepackage{latexsym}
\usepackage{amssymb}

\usepackage[pdftex]{graphicx}
\usepackage{amssymb}

\usepackage{url}

\usepackage{algorithm}
\usepackage{algorithmic}

\usepackage{xspace}

\usepackage{amsmath}
\usepackage{amsthm}

\theoremstyle{plain}
\newtheorem{theorem}{Theorem}[section]
\newtheorem{lemma}[theorem]{Lemma}

\newtheorem{corollary}[theorem]{Corollary}

\theoremstyle{definition}
\newtheorem{definition}[theorem]{Definition}
\newtheorem{construction}[theorem]{Construction}

\theoremstyle{definition}
\newtheorem{example}[theorem]{Example}

\newcommand{\mov}{\textnormal{Mov}}

\newcommand{\isomorphic}{\cong}

\newcommand{\darts}{\Gamma}

\newcommand{\opa}{\diamond} 
\newcommand{\opb}{\otimes}

\newcommand{\automorphismgroup}[1]{\textnormal{Aut}(#1)}

\newcommand{\objects}{\Omega}

\newcommand{\Rdisjoint}{\textnormal{(R1)}\xspace}
\newcommand{\RtripA}{\textnormal{(R2)}\xspace}
\newcommand{\RtripB}{\textnormal{(R3)}\xspace}

\newcommand{\Tidentity}{\textnormal{(T1)}\xspace}
\newcommand{\Tcycledisjoint}{\textnormal{(T2)}\xspace}
\newcommand{\Tfpf}{\textnormal{(T3)}\xspace}
\newcommand{\Ttransitive}{\textnormal{(T4)}\xspace}

\newcommand{\Ione}{\textnormal{(I1)}\xspace}
\newcommand{\Itwo}{\textnormal{(I2)}\xspace}

\newcommand{\fpf}{fixed-point-free\xspace}

\begin{document}
\title{\bf An enumeration of spherical \\ latin bitrades}
\author{\Large 
Ale\v s Dr\'apal\footnote{Supported by grant MSM 0021620839}\\
Department of Mathematics \\ Charles University \\ Sokolovsk\'a 83 \\ 186 75 Praha 8 \\ Czech Republic \\
~ \\
Carlo H\"{a}m\"{a}l\"{a}inen\footnote{Supported by Eduard \v Cech center, grant LC505.}  \\
Department of Mathematics \\ Charles University \\ Sokolovsk\'a 83 \\ 186 75 Praha 8 \\ Czech Republic \\
{\texttt carlo.hamalainen@gmail.com}\\
~ \\
Dan Rosendorf \\
Department of Mathematics \\
University of Wisconsin, Madison, USA
}
\date{}
\maketitle
\begin{abstract}
A {\em latin bitrade} $(T^{\opa},\, T^{\opb})$ is a pair of partial
latin squares which are disjoint, occupy the same set of non-empty
cells, and whose corresponding rows and columns contain the same set
of entries. A genus may be associated to a latin bitrade by constructing
an embedding of the underlying graph in an oriented surface. 
We report computational enumeration results on the number of spherical
(genus~0) latin bitrades up to size $24$.
\end{abstract}

\section{Introduction}

A {\em latin bitrade} $(T^{\opa},\, T^{\opb})$ is a pair of partial
latin squares which are disjoint, occupy the same set of non-empty
cells, and whose corresponding rows and columns contain the same set
of entries. One of the earliest studies of latin bitrades appeared 
in~\cite{DrKe1}, where they are referred to as {\em exchangeable partial
groupoids}. Latin bitrades are prominent in the study of 
{\em critical sets}, which are minimal defining sets of latin squares
(\cite{BaRe2},\cite{codose},\cite{Ke2}) and the intersections
between latin squares (\cite{Fu}). A genus may be associated to a latin
bitrade, and those of genus zero are known as {\em spherical} latin
bitrades. In this paper we report on the enumeration of the number of
$\tau$-isomorphism classes of spherical latin bitrades up to size $\left| T^{\opa} \right| = 24$.

For completeness we note that separated spherical latin
bitrades are equivalent to cubic $3$-connected bipartite planar
graphs~\cite{planareulerian}. Our enumeration method does not consider
equivalences where rows, columns, and symbols change their roles
(this corresponds to changing the $3$-colouring of the related cubic
3-connected bipartite graph), nor does it consider equivalences where
$(T^{\opa},\, T^{\opb})$ and $(T^{\opb},\, T^{\opa})$ are equivalent
(this corresponds to changing the $2$-colouring of the related cubic
3-connected bipartite graph). Thus our enumeration method counts a larger class of combinatorial objects compared to the
method in~\cite{holtonetc}.

\section{Latin bitrades}

A {\em partial latin square} $P$ of order $n > 0$ is an 
$n \times n$ array where each $e \in \{ 0, 1, \dots, n-1 \}$
appears at most once in each row, and at most once in each column.
A {\em latin square} $L$ of order $n > 0$ is an
$n \times n$ array where each $e \in \{ 0, 1, \dots, n-1 \}$
appears exactly once in each row, and exactly once in each column.
It is convenient to use setwise notation to refer to entries
of a (partial) latin square, and we
write $(i,j,k) \in P$ if and only if symbol $k$ appears in the
intersection of row $i$ and column $j$ of $P$.
In this manner, $P \subseteq A_1 \times A_2 \times A_3$ for finite sets
$A_i$, each of size $n$.
It is also convenient to interpret a (partial) latin square as a multiplication
table for a (partial) binary operator $\opa$, writing
$i \opa j = k$ if and only if $(i,j,k) \in T = T^{\opa}$.

\begin{definition}\label{defnBitradeA123}
Let $T^{\opa}$, $T^{\opb} \subseteq A_1 \times A_2 \times A_3$ be
two partial latin squares. Then $(T^{\opa},\, T^{\opb})$ is a {\em
bitrade} if the following three conditions are satisfied:
\begin{itemize}

\item[\Rdisjoint] $T^{\opa} \cap T^{\opb} = \emptyset$;

\item[\RtripA] for all $(a_1,\, a_2,\, a_3) \in T^{\opa}$ and all $r$,
$s \in \{1,\, 2,\, 3\}$,
$r \neq s$, there exists a unique $(b_1,\, b_2,\, b_3) \in T^{\opb}$
such that $a_r=b_r$ and $a_s=b_s$;

\item[\RtripB] for all $(a_1,\, a_2,\, a_3) \in T^{\opb}$ and all $r$,
$s \in \{1,\, 2,\, 3\}$,
$r \neq s$, there exists a unique $(b_1,\, b_2,\, b_3) \in T^{\opa}$
such that $a_r=b_r$ and $a_s=b_s$.

\end{itemize}
\end{definition}

Conditions~\RtripA and \RtripB imply that each row (column) of
$T^{\opa}$ contains the same subset of $A_3$ as the corresponding row
(column) of $T^{\opb}$.
A bitrade $(T^{\opa},\, T^{\opb})$ is {\em indecomposable\/}
if whenever $(U^{\opa},\,U^{\opb})$ is a bitrade such that
$U^{\opa}\subseteq
T^{\opa}$ 
and
$U^{\opb}\subseteq   
T^{\opb}$, then 
$(T^{\opa},\, T^{\opb})=  
(U^{\opa},\, U^{\opb})$.  
Bijections $A_i \rightarrow A'_i$, for $i = 1$, $2$, $3$, give an 
{\em isotopic} bitrade, and permuting each $A_i$ gives an {\em
autotopism}. We refer to the bijections
$A_1 \rightarrow A'_1$,
$A_2 \rightarrow A'_2$,
$A_3 \rightarrow A'_3$ as an {\em isotopism}.

In~\cite{Dr9,ales-geometrical} there is a representation of bitrades in terms of
three permutations $\tau_i$ acting on a finite set
(see also~\cite{hamalainen2007} for another proof).  
For $r \in \{1,\, 2,\, 3\}$, define the map 
$\beta_r \colon T^{\opb} \rightarrow
T^{\opa}$ where $(a_1,\, a_2,\, a_3) \beta_r = (b_1,\, b_2,\, b_3)$ if and only if
$a_r \neq b_r$
and $a_i = b_i$ for $i \neq r$.
By Definition~\ref{defnBitradeA123} each $\beta_r$ is a bijection.
Then
$\tau_1,\, \tau_2,\, \tau_3\colon T^{\opa} \rightarrow
T^{\opa}$ are defined by
\begin{align}
\tau_1 &= \beta_2^{-1}\beta_3, \qquad
\tau_2 = \beta_3^{-1}\beta_1, \qquad
\tau_3 = \beta_1^{-1}\beta_2. \label{eqnTau}
\end{align}
We refer to
$[\tau_1,\, \tau_2,\, \tau_3]$
as the $\tau_i$ {\em representation}.
We write $\mov(\pi)$ for the set of points that the (finite) permutation
$\pi$ acts on.

Generally we will assume that a bitrade is {\em separated}, that is,
each row, column, and symbol is in bijection with a cycle of $\tau_1$,
$\tau_2$, and $\tau_3$, respectively.

\begin{definition}\label{defnT1234}
Let $\tau_1$, $\tau_2$, $\tau_3$ be (finite) permutations
and let $\darts = \mov(\tau_1) \cup \mov(\tau_2) \cup \mov(\tau_3)$.  
Define four properties:
\begin{enumerate}

\item[\Tidentity] $\tau_1 \tau_2 \tau_3 = 1$;

\item[\Tcycledisjoint] if $\rho_i$ is a cycle of $\tau_i$
and $\rho_j$ is a cycle of $\tau_j$
then $\left| \mov(\rho_i) \cap \mov(\rho_j) \right| \leq 1$,
for any $1 \leq i < j \leq 3$;

\item[\Tfpf] each $\tau_i$ is \fpf;

\item[\Ttransitive] the group $\langle \tau_1,\, \tau_2,\, \tau_3 \rangle$ is
transitive on $\darts$.
\end{enumerate}

\end{definition}

By letting $A_i$ be the set of cycles of $\tau_i$, we
have the following theorem, which relates
Definition~\ref{defnBitradeA123}
and~\ref{defnT1234}.  

\begin{theorem}[{\cite{Dr9}}]\label{theoremDrapalTauStructure}
A separated bitrade $(T^{\opa},\, T^{\opb})$ is equivalent (up to isotopism) to three
permutations $\tau_1$, $\tau_2$, $\tau_3$ acting on a set $\darts$
satisfying \Tidentity, \Tcycledisjoint, and \Tfpf. 
If \Ttransitive is also satisfied then the bitrade is indecomposable.
\end{theorem}

To construct the $\tau_i$ representation for a bitrade we simply evaluate
Equation~\eqref{eqnTau}. In the reverse direction we have the following
construction:

\begin{construction}[$\tau_i$ to bitrade]\label{constructionTauToBitrade}
Let $\tau_1$, $\tau_2$, $\tau_3$ be permutations satisfying
Condition~\Tidentity, \Tcycledisjoint, and \Tfpf.
Let $\darts = \mov(\tau_1) \cup \mov(\tau_2) \cup \mov(\tau_3)$.  
Define $A_i = \{ \rho \mid \textnormal{$\rho$ is a cycle of $\tau_i$} \}$
for $i = 1$, $2$, $3$. Now define two arrays
$T^{\opa}$, $T^{\opb}$: 
\begin{align*}
T^{\opa}= \{( {\rho}_1,\, {\rho}_2,\, {\rho}_3 )\mid 
&\text{ $\rho_i \in A_i$ and 
$\left| \mov(\rho_1) \cap \mov(\rho_2) \cap \mov(\rho_3) \right| \geq
1$} \} \\
T^{\opb} = \{( {\rho}_1,\, {\rho}_2,\, {\rho}_3 ) \mid
&\text{ $\rho_i \in A_i$ and $x$, $x'$, $x''$ are distinct points of
$\darts$ such} \\
&\text{ that $x\rho_1=x'$, $x'\rho_2=x''$, $x''\rho_3=x$} \}. 
\end{align*}
By Theorem~\ref{theoremDrapalTauStructure} $(T^{\opa},\, T^{\opb})$ 
is a bitrade.
\end{construction}

\begin{example}\label{exampleIntercalateRep}
The smallest bitrade $(T^{\opa},\, T^{\opb})$ is the {\em intercalate},
which has four entries. The bitrade is shown below:
\begin{align*}
T^{\opa} = \begin{array}{c|cccc}
\opa & 0 & 1 \\
\hline 
0 & 0 & 1 \\
1 & 1 & 0
\end{array}
& \quad & 
T^{\opb} = \begin{array}{c|cccc}
\opb & 0 & 1 \\
\hline 
0 & 1 & 0 \\
1 & 0 & 1
\end{array}
\end{align*}
The $\tau_i$ representation is
$\tau_1 = (000,011)(101,110)$, $\tau_2 = (000,101)(011,110)$,
$\tau_3 = (000,110)(011,101)$, 
where we have written $ijk$ for $(i,j,k) \in T^{\opa}$ to make
the presentation of the $\tau_i$ permutations clearer.
By Construction~\ref{constructionTauToBitrade}
with $\darts = \{ 000, 011, 101, 110 \}$
we can convert the $\tau_i$ representation to a bitrade 
$(U^{\opa},\, U^{\opb})$:
\begin{align*}
U^{\opa} &= \begin{array}{c|cccc}
\opa & (000,101) & (011,110) \\
\hline 
(000,011) & (000,110) & (011,101) \\
(101,110) & (011,101) & (000,110)
\end{array} \\
U^{\opb} &= \begin{array}{c|cccc}
\opb & (000,101) & (011,110) \\
\hline 
(000,011) & (011,101) & (000,110) \\
(101,110) & (000,110) & (011,101)
\end{array}
\end{align*}
In this way we see that row $0$ of $T^{\opa}$ corresponds to 
row $(000,011)$ of $U^{\opa}$, which is the cycle
$( 000,011 )$ of $\tau_1$, and so on for the columns and symbols.
\end{example}
\begin{example}\label{exspherical}
The following bitrade is spherical:
\begin{align*}
T^{\opa} = 
\begin{array}{c|ccccc}
\opa & 0 & 1 & 2 & 3 & 4\\
\hline
 0 & 0 & ~ & 2 & ~ & 4\\ 1 & ~ & ~ & ~ & 4 & 2\\ 2 & 1 
& 3 & 0 & 2 & ~\\ 3 & 4 & 1 & ~ & 3 & ~\\ 4 & ~ & ~ & ~ & ~ & ~
\end{array}
& \quad & 
T^{\opb} = 
\begin{array}{c|ccccc}
\opb & 0 & 1 & 2 & 3 & 4\\
 \hline 0 & 4 & ~ & 0 & ~ & 2\\ 1 & ~ & ~ & ~ & 2 & 4\\ 2 & 0 & 1 & 2 & 3 & ~\\ 3 & 1 & 3 & ~ & 4 & ~\\ 4 & ~ & ~ & ~ & ~ & ~\\\end{array}
\end{align*}
Here, the $\tau_i$ representation is
\begin{align*}
\tau_1 &= (000, 022, 044)(134, 142) (201,213,232,220) (304,333,311) \\
\tau_2 &= (000,304,201)(213,311)(022,220)(134,232,333)(044,142) \\
\tau_3 &= (000,220)(201,311)(022,232,142)(213,333)(044,134,304)
\end{align*}
\end{example}

If $Z = [\tau_1,\, \tau_2,\, \tau_3]$ then the inverse is denoted by
$Z^{-1} = [\tau_1^{-1},\, \tau_2^{-1},\, \tau_2 \tau_1]$. This is
equivalent to exchanging $(T^{\opa},\, T^{\opb})$ for $(T^{\opb},\,
T^{\opa})$ and relabelling entries as the following lemma shows:

\begin{lemma}
Let $(T^{\opa},\, T^{\opb})$ be a separated bitrade with representation
$Z = [\tau_1,\, \tau_2,\, \tau_3]$. Then the inverse bitrade
$(T^{\opb},\, T^{\opa})$ has representation, denoted $Z^{-1}$. This representation is isomorphic to
$[\tau_1^{-1},\, \tau_2^{-1},\, \tau_2 \tau_1]$.  
\end{lemma}

\begin{proof}
Observe that Definition~\ref{defnBitradeA123} does not specify the
order of the two partial latin squares in a bitrade so the bitrade
$(T^{\opb},\, T^{\opa})$ is well-defined.  
For this bitrade we have permutation representation
$[\nu_1,\, \nu_2,\, \nu_3]$. Using Equation~\eqref{eqnTau} 
we find that
$\nu_1 = \beta_2 \beta_3^{-1}$
and
$\nu_2 = \beta_3 \beta_1^{-1}$. Taking the conjugate of $\nu_r$ by 
$\beta_3^{-1}$ gives us a representation that
acts on the set $T^{\opa}$. In particular,
$\beta_3^{-1} \nu_1 \beta_3 = \beta_3^{-1} \beta_2 \beta_3^{-1} \beta_3 = 
\beta_3^{-1} \beta_2 = 
(\beta_2^{-1} \beta_3)^{-1} = \tau_1^{-1}$
and
$\beta_3^{-1} \nu_2 \beta_3 = \beta_3^{-1} \beta_3 \beta_1^{-1} \beta_3 = 
\beta_1^{-1} \beta_3 = 
(\beta_3 \beta_1^{-1})^{-1} = \tau_2^{-1}$. This shows that the representation
$[\nu_1,\, \nu_2,\, \nu_3]$
is isomorphic to
$[\tau_1^{-1},\, \tau_2^{-1},\, \tau_2 \tau_1]$.
\end{proof}

\begin{figure}
\begin{center}
\includegraphics{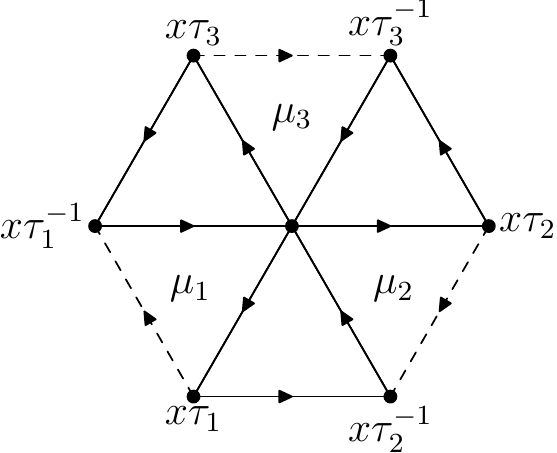}
\end{center}
\caption{Rotation scheme for the combinatorial surface. Dashed
lines represent one or more edges, where $\mu_i$ is a cycle of $\tau_i$.
The vertex in the centre is $x$.}
\label{figDrapalRotationScheme}
\end{figure}

Let $D$ be the graph with vertices $\darts$ and directed edges
$(x, y)$ where $x \tau_i = y$ for some $i$. A rotation scheme can be
imposed on this graph, by ordering the edges as shown in
Figure~\ref{figDrapalRotationScheme}. This
turns $D$ into an oriented combinatorial surface.
Define
$\textnormal{order}(\tau_i)  = z(\tau_1) + z(\tau_2) + z(\tau_3)$, the total number
of cycles, and
$\textnormal{size}(\tau_i) = \left| \darts \right|$. 
By some basic counting arguments we find that there are 
$\text{size}(\tau_i)$ vertices, $3 \cdot \text{size}(\tau_i)$ edges, and 
$\text{order}(\tau_i) + \text{size}(\tau_i)$ faces. Then Euler's formula
$V-E+F=2-2g$ gives
\begin{equation}\label{eqnDrapalGenus}
\text{order}(\tau_i) = \text{size}(\tau_i) + 2 - 2g
\end{equation}
where $g$ is the genus of the combinatorial surface. We say that
the bitrade $T$ has genus $g$.
If $g = 0$ then we say that the bitrade is {\em spherical}.

\section{Slide expansion}

Let $[\tau_1,\, \tau_2,\, \tau_3]$ be a spherical bitrade on the
set $\darts$. Choose some $x \in \darts$ and fix a direction $j \in
\{1,2,3\}$. Set $k = j+1$ and $l = k+1$ (calculating modulo~$3$).
Suppose that (1) the $j$-cycle at $x$ has length at least~$3$
and (2)
the $k$-cycle through $x$ and the $l$-cycle through
$w = x \tau_j$
do not meet in any common point. Then the
{\em slide expansion} at $x \in \darts$ in direction $j$ is the following
augmentation of the cycles of $\tau_1$, $\tau_2$, and $\tau_3$ by a new
point $u \notin \darts$:
\begin{align*}
\tau_j &\colon (a, x, w, b, \dots) \mapsto (a, u, b, \dots) \textnormal{ and } (x, w) \\
\tau_k &\colon (\dots, x, z, \dots) \mapsto (\dots, x, u, z, \dots) \\
\tau_l &\colon (\dots, y, w, \dots) \mapsto (\dots, y, u, w, \dots).
\end{align*}
We denote the result of the slide expansion on $X = [\tau_1,\, \tau_2,\, \tau_3]$ at vertex
$x$ in direction $j$ by $\textnormal{slide}(X, x, j)$.

A separated bitrade $[\tau_1,\, \tau_2,\, \tau_3]$ is {\em bicyclic} if
there exists $j$ such that $\tau_j$ consists of only two cycles. These
bicyclic bitrades form the root nodes of our enumeration procedure due
to the following result:

\begin{theorem}[\cite{ales-spherical}]
Let $Z = [\tau_1,\, \tau_2,\, \tau_3]$ be a spherical bitrade. Then 
there exists a sequence $W_0$, $W_1$, \dots, $W_{\ell} = Z$ where
$W_0$ is a bicyclic bitrade and
$W_{i+1}$ is the result of a slide expansion on
either $W_{i}$ or $W^{-1}_i$, for $0 \leq i < \ell$.
\end{theorem}

The opposite of a slide expansion is a 
{\em slide contraction}. For a spherical bitrade 
$[\tau_1,\, \tau_2,\, \tau_3]$ on the
set $\darts$, choose some $u \in \darts$ and fix a direction $j \in
\{1,2,3\}$. Set $k = j+1$ and $l = k+1$ (calculating modulo~$3$).
The slide contraction at $u$ in direction $j$ modifies the cycles of
$\tau_1$, $\tau_2$, and $\tau_3$ as follows:
\begin{align*}
\tau_j &\colon (a, u, b, \dots) \textnormal{ and } (x, w) \mapsto (a, x, w, b, \dots) \\
\tau_k &\colon (\dots, x, u, z, \dots) \mapsto (\dots, x, z, \dots) \\
\tau_l &\colon  (\dots, y, u, w, \dots) \mapsto (\dots, y, w, \dots).
\end{align*}
A necessary condition for a slide contraction to be valid is that the
$u$-cycle in $\tau_k$ and $\tau_l$ has length at least $3$ (this avoids
the creation of a fixed point when $u$ is removed). In general it is
necessary to check the bitrade conditions \Tidentity, \Tcycledisjoint,
and \Tfpf to ensure that a slide contraction is valid.

\section{Enumeration}

To enumerate separated spherical bitrades we must fix the notion of $\tau$-isomorphism.
Two bitrades $[\tau_1,\, \tau_2,\, \tau_3]$ and $[\nu_1,\, \nu_2,\, \nu_3]$
on the same set $\darts$ are $\tau$-isomorphic if there is a permutation
$\theta \in \textnormal{Sym}(\darts)$ such that
$\tau_i^{\theta} = \nu_i$ for all $i \in \{1,2,3\}$. A $\tau$-automorphism of
a bitrade is defined similarly.

We follow the presentation in~\cite{classificationalgorithms}.
In the context of an algorithm for enumerating combinatorial objects,
the {\em domain} of a search is the finite set $\objects$ that contains
all objects considered by the search. The search space is conveniently
modelled by a rooted tree, with nodes corresponding to elements of
$\objects$, and nodes joined by an edge if they are related by one
search step. The root node is the starting node of the search. We write 
$C(X)$ for the set of all child nodes of $X$, and $p(X)$ for the parent
node of $X$.

To reduce the search time we use the method of {\em canonical
augmentation}~\cite{mckay}, following the presentation 
in~\cite{classificationalgorithms}. In terms of the search tree, the ordered
pair $(X, p(X))$ characterises the augmentation used to generate $X$
from the node $p(X)$ during the search. In our case, this augmentation
$(X, p(X))$ is the result of a slide expansion on $p(X)$.  
The goal of the canonical augmentation algorithm is to choose one of the
possible parents of $X$ to be the canonical parent. Formally,
for any nonroot node $X$,
the {\em canonical parent} $m(X) \in \objects$
must satisfy: 
\begin{enumerate}

\item[(C1)] for all nonroot objects $X$, $Y$ it holds that
$X \isomorphic Y$ implies 
$(X, m(X)) \isomorphic (Y, m(Y))$.

\end{enumerate}
The next property captures the canonical augmentations in
the search tree:

\begin{enumerate}

\item[(C2)] a node $Z$ occurring in the search tree is generated by
canonical augmentation if $(Z, p(Z)) \isomorphic (Z, m(Z))$.

\end{enumerate}
Algorithm~\ref{alg1} gives pseudo-code for the general canonical
augmentation search algorithm.  Algorithm~\ref{alg2} is a more explicit
version for enumerating separated spherical latin bitrades.
The correctness of Algorithm~\ref{alg1} relies on two properties:
\begin{enumerate}

\item[\Ione] for all nodes $X$, $Y$ it holds that if $X \isomorphic Y$,
then for every $Z \in C(X)$ there exists a $W \in C(Y)$ such that
$(Z,X) \isomorphic (W,Y)$.

\item[\Itwo] for every nonroot node $X$, there exists a node $Y$ such
that $X \isomorphic Y$ and $(Y,m(Y)) \isomorphic (Y,p(Y))$.
\end{enumerate}

Condition \Ione says that isomorphic nodes have isomorphic children.
Condition \Itwo says that every nonroot node is generated by some
canonical augmentation.  

\begin{algorithm}
\caption{CANAUG-TRAVERSE($X$)}
\label{alg1}
\begin{algorithmic}[1]
\STATE report $X$ (if applicable)
\FORALL{$\mathcal{Z} \in \{ C(X) \cap \{ aZ \colon a \in \automorphismgroup{X} \} \colon Z \in C(X) \}$}
    \STATE select any $Z \in \mathcal{Z}$
    \IF{$(Z,p(Z)) \isomorphic (Z,m(Z))$}
        \STATE CANAUG-TRAVERSE($Z$)
    \ENDIF
\ENDFOR
\end{algorithmic}
\end{algorithm}

\begin{algorithm}
\caption{CANAUG-SPHERICAL($X$)}
\label{alg2}
\begin{algorithmic}[1]
\STATE report $X$ (if applicable)
\FORALL{$\mathcal{Z} \in \{ C(X) \cap \{ aZ \colon a \in \automorphismgroup{X} \} \colon Z \in C(X) \}$}
    \STATE select any $Z \in \mathcal{Z}$
    \IF{$(Z,p(Z)) \isomorphic (Z,m(Z))$}
        \STATE CANAUG-SPHERICAL($Z$)
        \IF{$Z^{-1}$ has no parents}
            \STATE CANAUG-SPHERICAL($Z^{-1}$)
        \ENDIF
    \ENDIF
\ENDFOR
\IF{$(\left| X \right| + 1) \equiv 0 \pmod{2}$}
    \FORALL{bicyclic $Z$ with $\left| Z \right| = \left| X \right| + 1$}
        \STATE CANAUG-SPHERICAL($Z$)
    \ENDFOR
\ENDIF
\end{algorithmic}
\end{algorithm}

\begin{theorem}[\cite{classificationalgorithms}]
When implemented on a search tree satisfying \Ione and \Itwo, the algorithm
\textnormal{CANAUG-TRAVERSE} reports exactly one node from every isomorphism class of
nodes.  
\end{theorem}

\subsection{Canonical form}

The canonical form $\hat{Z}$ is a relabelling of $Z$
such that $\hat{Z} = \hat{Y}$ if and only if $Z \isomorphic Y$. We need
to be able to efficiently compute the canonical form in order to
construct the canonical parent of a bitrade.

To compute the canonical form of a bitrade $Z = [\tau_1,\, \tau_2,\,
\tau_3]$, we fix a starting point $x \in \darts$. We then perform a
breadth-first traversal of the bitrade's underlying directed graph.
Vertices are labelled from the set $\{ 1, \dots, \left| Z \right| \}$ in
the order that they are visited. At a vertex $v$ we first visit the 
$\tau_1$ cycle at $v$, then the
$\tau_2$ cycle at $v$, and then the 
$\tau_3$ cycle at $v$ (only if the cycle has not been visited before).
As we traverse a cycle we label unseen vertices and append the cycle
(using the new vertex labels) to
an array $C$. After traversing a cycle we append a marker (here, $-1$)
to denote the end of the cycle. While traversing a 
$\tau_i$ cycle at the vertex $v$, we check if the neighbouring 
$\tau_{j}$ and $\tau_{k}$ cycle at $v$ has been visited ($j$, $k \neq i$).
If it has not, we push $(j,v)$ and $(k,v)$ onto the queue, making sure
that $j < k$. The lexicographically least $C = C(x)$ is
the canonical form $\hat{Z}$ of $Z$.
Since $\hat{Z}$ encodes just the
cycle structure of the $\tau_i$ permutations, it follows that
$\hat{Z} = \hat{Z'}$ if and only if $Z \isomorphic Z'$.

The canonical parent $m(Z)$ is now computed as follows. First we find
the canonical form $\hat{Z}$ of $Z$. In our case, we choose the
lexicographically maximal pair $(j, u)$ such that a slide contraction
may be performed on $\hat{Z}$ at $u$ in direction $j$. Let $Z_0$ be the
result of this contraction. Then $m(Z) = \kappa^{-1}(Z_0)$ where
$\kappa \colon Z \mapsto \hat{Z}$ is the canonical relabelling map.

\subsection{Correctness}

Condition \Ione holds since a slide expansion 
$\textnormal{slide}(X, x, j)$ is equivalent to 
$\textnormal{slide}(Y, y, j)$ where $X \isomorphic Y$ and $y$ is the
image of $x$ under the isomorphism. In other words, the slide expansion
relies only on the cycle structure of the bitrade, not its specific
labelling. Condition \Itwo holds because the canonical parent $m(Z)$ is chosen
only based on the canonical form of $Z$ and so is independent of
isomorphism.

Algorithm~\ref{alg2} produces bicyclic bitrades at each step where
even-sized bitrades are being generated. The following lemma
and corollary are necessary to show that no $\tau$-isomorphism class of
bicyclic bitrades is produced more than once.

\begin{lemma}[Lemma~10.1~\cite{ales-spherical}]\label{bicyclic_structure}
Let $[\tau_1,\, \tau_2,\, \tau_3]$ be a bitrade such that
$\tau_j$, for some $j \in \{1,2,3\}$, consists of exactly two cycles.
Let one of the cycles be $(x_0, \dots, x_{n-1})$. Then the other cycle
of $\tau_j$ can be expressed as
$(x'_{n-1}, \dots, x'_0)$ in such a way that
$\tau_{j+1}$ consists of cycles $(x_i, x'_i)$ and $\tau_{j-1}$ consists
of cycles $(x_{i+1}, x'_i)$, for $0 \leq i < n$.
\end{lemma}

\begin{corollary}\label{corbicyclic}
Let $Z$ be a bicyclic bitrade. Then there is no $X$ such that
$\textnormal{slide}(X, x, j) = Z$ or
$\textnormal{slide}(X, x, j) = Z^{-1}$ for any $x$, $j$.
\end{corollary}

\begin{proof}
The slide contraction shortens the cycle at some point $u$
in $\tau_{j+1}$ and $\tau_{j-1}$.
By Lemma~\ref{bicyclic_structure} at least one of the cycles in
$\tau_{j+1}$ or $\tau_{j-1}$ has length $2$, so the slide contraction
would introduce a fixed point, contradicting \Tfpf.
\end{proof}

It can be shown that for any bicyclic bitrade $Z$ we have $Z \isomorphic
Z^{-1}$.  Corollary~\ref{corbicyclic} shows that lines 11--15 of
Algorithm~\ref{alg2} are correct.

\section{Implementation and relation to other work}

See \cite{code} for C++ computer code. The MPI library is used to allow
parallel processing of the search tree on a Linux cluster.
The size of the automorphism group $\automorphismgroup{Z}$ is bounded
by the size of $\darts$. Further, the number of homogeneous bitrades
is quite small, so usually the automorphism group is even smaller.
We just compute the elements of $\automorphismgroup{Z}$ explicitly.

Cavenagh and Lisonek~\cite{planareulerian} showed that spherical
bitrades are equivalent to planar Eulerian triangulations. To verify our
results (for small $n$) we used plantri~\cite{plantri} to generate
planar Eulerian triangulations. We then applied a simultaneous vertex
$3$-colouring and face $2$-colouring, giving a bitrade $(T^{\opa},\,
T^{\opb})$. Then each conjugate of $(T^{\opa},\, T^{\opb})$ was
produced, and finally isomorphic copies removed.

The number of $\tau$-isomorphism classes of spherical bitrades with 
$\left| T^{\opa} \right| = n$ is given below:

\begin{center}
\begin{tabular}{|c|c|c|}
\hline $n$ & \# iso classes & $e_3(n_v)$ \\
\hline 
\hline 4 & 1 & 1 \\
\hline 6 & 3 & 1 \\
\hline 7 & 1 & 1 \\
\hline 8 & 6 & 2 \\
\hline 9 & 9 & 2 \\
\hline 10 & 30 & 8 \\
\hline 11 & 51 & 8 \\
\hline 12 & 198 & 32 \\
\hline 13 & 470 & 57 \\
\hline 14 & 1623 & 185 \\
\hline 15 & 4830 & 466 \\
\hline 16 & 16070 & 1543 \\
\hline 17 & 51948 & 4583 \\
\hline 18 & 175047 & 15374 \\
\hline 19 & 588120 & 50116 \\
\hline 20 & 2015226 & 171168 \\
\hline 21 & 6933048 & 582603 \\
\hline 22 & 24123941 & 2024119 \\
\hline 23 & 84428820 & 7057472 \\
\hline 24 & 297753519 & 24873248 \\
\hline
\end{tabular}
\end{center}

The third column is the value $e_3(n_v)$, which is 
the number of isomorphism classes of eulerian plane triangulations with
connectivity at least $3$ (see Table~6 in \cite{fastgen} for values
$n_f = 2n$ in their first column and also~\cite{A000108}). The enumeration resulting in 
$e_3(n_v)$ considers isomorphic any two graphs that are related by a
permutation of the face 2-colouring, changing the orientation of the
permutations, and permuting the cyclic ordering. So for large values of
$n$ there will be at most $2 \cdot 2 \cdot 3 = 12$ ways to relabel one
of their graphs to obtain one of our bitrades. For example, with $n =
24$ we have $e_3(n_v) = 24873248$ and 
$12 \cdot 24873248 = 298478976 \approx 297753519$.
Also, Wanless~\cite{wanlessenumeration} enumerated spherical bitrades
up to size $19$, under various various equivalences, but not
$\tau$-isomorphism as studied in this paper.


\end{document}